\documentclass[12pt,reqno]{amsart}
\usepackage{latexsym}
\usepackage{amssymb}

\usepackage{epsfig}

\renewcommand{\epsilon}{\varepsilon}

\newcommand{\PP}{{\mathbb P}}

\newcommand{\C}{{\mathbb C}}

\newcommand{\CP}{\C\PP}

\newcommand{\E}{{\mathbf E}}

\newcommand{\half}{{\frac{1}{2}}}

\renewcommand{\phi}{\varphi}

\newcommand{\ecal}{\mathcal{E}}

\newtheorem{theo}{{\sc Theorem}}

\newtheorem{lem}[theo]{{\sc Lemma}}

\title[Addendum to   ``Energies of zeros of random sections on Riemann surfaces''
 ] {Addendum to  ``Energies of zeros of random sections on Riemann surfaces''. Indiana Univ. Math. J. 57 (2008), no. 4, 1753--1780
 }

\author{Qi Zhong}

\author{Steve Zelditch}
\address{Department of Mathematics, Northwestern University\\
Evanston, IL, USA} \email{zelditch@math.northwestern.edu}

\date{\today}

\begin{document}

\maketitle

The purpose of this note is to resolve an apparent
 discrepancy between the calculations in the article \cite{ABS} of  Armentano- Beltran-Shub    (henceforth ABS) and that in
  Qi Zhong's article \cite{Zh}
of the asymptotics of the expected  energy
of zeros of random polynomials on the Riemann sphere $S^2$ with respect to the log chordal distance $\log [z, w]$. For the sake of brevity, we
do not repeat the statement of the problem but refer to \cite{Zh,ABS} for the background and notation.
We show that the calculation in \cite{Zh} which uses
  the general  method of Green's function and correlation functions  gives the same answer as in \cite{ABS}.

The discrepancy between \cite{ABS} and \cite{Zh}  arose because \cite{Zh} actually contains two calculations of the asymptotic energy with respect to
$\log [z, w]$, one explicit and one implicit. The implicit one is indicated in the Remarks after \cite{Zh} Theorem 1.2.
There, it is pointed out that  the general  Green's energy method of the paper
applies to the log chordal energy on $S^2$  if one adds a constant to  $\log [z,w]$ to convert it to the Green's function.
 The Green's energy asymptotics in \cite{Zh} are correct and their application to the log chordal energy does give
 the correct answer, as we verify in this note.  But  the energy asymptotics were not computed by that method in \cite{Zh}.
Rather in Section 5.2 of \cite{Zh}, the energy asymptotics were calculated without converting $\log [z, w]$ to the Green's function,
and there are some errors in the calculation of the integrals that arise.  As a result, the  asymptotics stated in   Theorem 1.3 (2) (see (1.13))
are incorrect. The correct asymptotics are given here.

Besides correcting the calculation of the $\log [z,w]$ energy in \cite{Zh}, the purpose of this note is to clarify
what is correct and what is incorrect in  \cite{Zh}. Most importantly, the main result of   \cite{Zh} is correct,
and proves the correct  asymptotics for the  Green's energy of
zeros of random holomorphic sections of powers of positive Hermitian line bundles for
all Hermitian metrics of positive $(1,1)$ curvature over all Riemann surfaces.  By comparison, the later result of \cite{ABS}
 only give the asymptotics in one special case, the round metric on $S^2$. However, as pointed out to us by M. Shub, the formula in
 \cite{ABS} is exact: the $o(N)$ error term in Zhong's formula is zero. We do not prove that here, because we derive the asymptotics
 from a general formula for metrics on Riemann surfaces, where in general the $o(N)$ term is non-zero.

The discrepancy between the asymptotics stated in Theorem 1.3 (2) of \cite{Zh}  and those of \cite{ABS} was spotted
by D. Hardin and E. Saff while Q. Zhong was a postdoc at Vanderbilt. They also informed the authors of \cite{ABS} about the
apparent discrepancy. A reference is made in \cite{ABS} to an erratum by Zhong. The present note supplants his erratum.

\subsection{ABS versus Zhong}

The calculation in question concerns the expected  logarithmic energy  of zero sets
of polynomials of degree $N$  on the Riemann sphere. The first important point of comparison is the normalizations
of the energy in \cite{ABS} and \cite{Zh}.

In both articles, the Riemann sphere is identified with  the sphere $S(\half) $ of radius $\half$ centered at
$(0, 0, \half)$. The logarithm chordal energy is then given by
$$ - \log [z,w], \;\;\; [z,w] = \frac{|z - w|}{\sqrt{1 + |z|^2} \sqrt{1 + |w|^2}}. $$

 For this kernel on $S(\half)$,
the  expected log  energy calculated in \cite{ABS} is
\begin{equation} \label{SHUB} \mbox{ABS}\;\;\;\E \ecal_{ABS}^N \sim \frac{N^2}{4}  - \frac{N \log N}{4} - \frac{N}{4}. \end{equation}


 Qi Zhong calculated the energy not for $\log [z,w]$ but for minus the Green's function of $S(\half)$. Following \cite{Zh},
 we denote minus the Green's function of a metric $g$  by $G_g$ and denote the usual Green's function by $\hat{G}_g$.
 Thus, $G_g = - \hat{G}_g$. We denote the Green's function of $S(\half)$ by $\hat{G}_{\half}$.  Then  $-G_{\half}(z,w): = - \frac{1}{2 \pi} \log [z,w] - C_{\half}$
 for a certain constant $C_{\half}$. For the  $G_g$-Green's energy, Zhong prove the general (and correct) asymptotics,
\begin{equation} \label{ZH} \mbox{ZHONG} \;\;\;\;\E \ecal^N_{ZH} = - \frac{1}{4 \pi} N \log N - \frac{N}{4 \pi} - N
\int F_g(z,z) \omega_h/\pi + o(N), \end{equation}
where $F_g$ is the Robin constant of $g$. We denote the Robin constant of $S(\half)$ by $F_{\half}$;
 note that it is a constant for the round metric.

To compare the two formulae (\ref{SHUB})-(\ref{ZH}), we need to calculate the constants $C_{\half}$ and $F_{\half}$,
add $C_{\half}$ to $ G_{\half}$ to convert it to $- \frac{1}{2 \pi} \log [z,w]$, and substitute the value of $F_{\half}$. In addition, as detailed in \S \ref{NORM},
 we need to multiply
 by $\frac{1}{\pi}$ to convert ABS to ZHONG. Thus, the following Lemma asserts that the calculations of \cite{Zh} and \cite{ABS} agree:

 \begin{lem} \label{CONSTANTS} $C_{\half} = \frac{1}{4 \pi}$ and  $F_g(z,z) = - \frac{1}{4 \pi}$. Hence,
 $$\pi \left(\mbox{ZHONG} + C_{\half} N(N - 1) \right) = \mbox{ABS}. $$
 \end{lem}

\subsection{\label{NORM} Normalizations}

To make the normalizations in \cite{ABS} and \cite{Zh} consistent we need to observe that:

\begin{itemize}

\item  (i) Zhong defines the energy using  $- \frac{1}{2 \pi} \log [z,w] - C_{\half}$
while ABS use $- \log [z,w]$.

\item (ii)
ABS sums over $i < j$ while Zhong sums over $i \not= j$.

\end{itemize}

 Zhong assumes that the Riemannian area form is $dV = \omega_h$ with $\int_{\CP^1} \omega_h = \pi$;  see (2.2) of \cite{Zh}. Since this
is  the area of  the sphere of radius $\half$,  \cite{Zh} and
ABS are working on $S(\half)$.

To emphasize that the  metric quantities pertain to the sphere of radius  $\half$, we subscript
all metric quantities by the respective radius, except for the geodesic distance, which we denote by $r$.


\subsection{Proof of Lemma \ref{CONSTANTS}}

 Both constants  depend on the radius we pick for $S^2$, namely radius $\half$. To  keep track of constants,
 we denote by $A$ the area of $S^2$ in the given  metric. For the round metric of area $\pi$,  $G_{\half}(z,w) $ is a function of
the geodesic distance $r(z,w)$ and hence of the chordal distance $[z,w]$. Subscripting with the radius, we have
$$\hat{G}_{\half}(z,w) = \frac{1}{2 \pi} \log [z,w]_{\half} + C_{\half}.$$

\subsection{The constant $C_{\half}$}

  The constant is
determined by the fact that $\int_{\CP^1} G_{\half}(z,w) \omega_w = 0$ for
all $z$, which becomes
$$\begin{array}{lll} A(S(\half)) C_{\half} &= & - \frac{1}{2 \pi} \int_0^{2 \pi} \int_0^{\frac{\pi}{2}} \log [0, r]_{\half}
(\half \sin 2 r dr) d \theta
= - \int_0^{\frac{\pi}{2}} \log [0, r]_{\half} (\half \sin 2 r dr)  \\ && \\
& = & - \int_0^{\pi} \log (\half \sin \phi) (\frac{1}{4} (\sin \phi) d\phi)\\ && \\
& = & - (\frac{1}{4} \int_0^{\pi} \log ( \sin \phi)  (\sin \phi) d\phi)  - (\frac{1}{4} (\log \half) \int_0^{\pi}  (\sin \phi) d\phi) .
 \\ && \\
& = & - (\frac{1}{4} \int_0^{\pi} (\log \sqrt{2(1 - \cos  \phi)})  (\sin \phi) d\phi)  - (\frac{1}{4} (\log \half) \int_0^{\pi}  (\sin \phi) d\phi)
 \\ && \\
& = & - \frac{1}{4} \log \frac{4}{e}  - (\frac{1}{4} (\log \half) \int_0^{\pi}
(\sin \phi) d\phi) =  - \frac{1}{4} \log \frac{4}{e}  - (\frac{1}{4} (\log \half) 2.

 \end{array}$$
 We conclude that
 \begin{equation} C_{\half} =  - \frac{1}{4  \pi} \log \frac{4}{e}  - (\frac{1}{4\pi} (\log \half) 2 = \frac{1}{4 \pi}. \end{equation}

 Then Zhong's minus Green's function is given by,
 \begin{equation}  G_{\half}(z,w) = -  \frac{1}{2 \pi} \log [z,w] +  \frac{1}{4  \pi} \log \frac{4}{e}  + (\frac{1}{4\pi} (\log \half) 2
 = -  \frac{1}{2 \pi} \log [z,w] - \frac{1}{4 \pi}. \end{equation}
 It follows that,
 \begin{equation} \label{FOLLOWS}  -  \frac{1}{2 \pi} \log [z,w] =  G_{\half}(z,w)  + \frac{1}{4 \pi}. \end{equation}

\subsection{Robin constant}

We further show that
\begin{equation} \label{FG} F_{\half}(z,z) = - \frac{1}{4 \pi}. \end{equation}
Here,  $F_{\half}(z,z) $ is  the constant in  the expansion
$$G_{\half}(z,w) = - \frac{1}{2 \pi} \log r + F_{\half}(z,z) + O(r), \;\;\; r \to 0. $$
In fact, it is the same constant that we just calculated.


Indeed,  $\hat{G}_{\half} = \frac{1}{2 \pi} \log (\half \sin 2 r) +\frac{1}{4 \pi}, $ and $\half \sin 2 r  = r f(r), f(0) = 1$. So $\log f(0) = 0$
and $$F_{\half} = - \frac{1}{4 \pi}, $$as desired.

\subsection{Conclusion}

Adding $\frac{1}{4 \pi}$ to $G_{\half}$  in (\ref{FOLLOWS}) results in adding
 $\frac{1}{4 \pi} N(N - 1)$ to Zhong's asymptotics (\ref{ZH}). We then substitute  (\ref{FG}) to find that (\ref{ZH}) equals
$$ \frac{1}{4 \pi} N(N - 1)  - \frac{1}{4 \pi} N \log N - \frac{N}{4 \pi} - N
(- \frac{1}{4 \pi}) + o(N). $$

Finally, we  multiply by $\pi$, to get
$$\frac{N^2}{4}  - \frac{N \log N}{4} - \frac{N}{4} + o(N), $$
which is the same as (\ref{SHUB}). This
completes the proof of Lemma \ref{CONSTANTS}, and proves that
the calculations of  \cite{ABS} and \cite{Zh} agree.

Finally, we thank B. Shiffman for checking over the calculations.

\end{document}